\magnification=1200
\hsize =15true cm
\baselineskip = 12pt
\vsize= 22true cm

\centerline{\bf $\Gamma$-convergence of Integral  Functionals}

\centerline{\bf Depending on Vector-valued Functions over
Parabolic Domains}

 \vskip 0.3cm
 \centerline{Huai-Yu Jian}

\centerline{( Department of  Mathematics, Tsinghua University,
Beijing 100084, P.R. China )}

\vskip 0.4cm

{\bf ABSTRACT:} We study $\Gamma$-convergence for a sequence of
parabolic functionals, $F^\varepsilon (u)=\int_0^T\int_{\Omega}
f({x\over\varepsilon}, t, \nabla u) dxdt$ as $\varepsilon \to 0$,
where the integrand $f$ is nonconvex, and periodic on the first
variable. We obtain the representation formula of the
$\Gamma$-limit. Our results in this paper support a conclusion
which relates  $\Gamma$-convergence of parabolic functionals to
the associated gradient flows and confirms one of De Giorgi's
conjectures partially.

\vskip 0.3cm


{\bf 1991 MR Classification No.}: 35B27,  \ \ 49J45

\

\
\centerline { \bf 1.   INTRODUCTION }

We begin with the characterization of $\Gamma$-convergence in [1, 2].

{\bf  DEFINITION 1.1.}   {\sl Let $(X,\tau )$ be a first countable topological
space and $\lbrace F^h \rbrace_{h=1}^{\infty} $ be a sequence of functionals
from $X$
to $\bar R = R\cup \lbrace - \infty, \infty \rbrace $,\ \
$u \in X, \lambda \in
\bar R .$  We call $$ \lambda = \Gamma (\tau ) \lim_{h \to \infty }F^h (u)$$
if and only if for every sequence $\{u^h\}$ converging to u in $(X,\tau )$
$$ \lambda \le \liminf_{h  \to \infty }F^h(u^h) ,  \eqno (1.1)$$
\noindent and there exists a sequence  $\{u^h\} $  converging to
$ u$ in $(X,\tau )$
such that
$$ \lambda \ge \limsup_{h  \to \infty }F^h(u^h) .  \eqno (1.2) $$
\noindent We call $ \lambda = \Gamma (\tau ) \lim_{\varepsilon \to a }F^\varepsilon
(u) $ if and only if for every $ \varepsilon _h \to  a $ \ \ $(h \to \infty )$
$$ \lambda = \Gamma (\tau )\lim_{h \to \infty }F^{\varepsilon _h}(u).$$}

Throughout this paper, we assume that
 $\Omega$ is  a bounded open set  in $ R^n$. Let  $p>1$, $T>0$, and $m$ be
 a positive integer. Denote
 $$\Omega _T=\Omega\times (0,T), \ \
V_p(\Omega _T,m)=L^P([0 , T],W^{1,p}(\Omega , R^m)), $$
$$V_p^0(\Omega _T ,m)=L^p([0,T],W_0^{1,p}(\Omega , R^m))  ,$$ and
$$ Du(x,t)=\nabla u (x, t)=\bigl ({\partial u^i(x,t) \over \partial x_j } \Bigr ) \quad
(1 \leq i \leq m , 1 \leq j \leq n )  $$ for a vector valued function $u$.

 Consider the fuctionals
$$F_1^\varepsilon (v, \Omega ) =\int _{\Omega }f_1 ({x\over \varepsilon}
,Dv)dx , \ \ v\in W^{1,p}(\Omega , R^m) , (\varepsilon \to 0^+) \eqno (1.3)$$
and the corresponding parabolic
 functionals in the following form:
$$F^\varepsilon (u, \Omega _T) =\int _{\Omega _T}f ({x\over \varepsilon},t,
,Du)dxdt ,\ \ u\in V_p (\Omega _T, m),
(\varepsilon \rightarrow 0^+ ),  \eqno (1.4) $$
where  $f\colon R^{n+1} \times R^{mn} \rightarrow R $ is a
{\sl Caratheodory } function satsfying
$$ C_1 |\lambda |^p \le f (x,t,\lambda ) \le C_2 (1+ |\lambda |^p )
 \eqno (1.5)$$
 for some positive constants $C_2 >C_1$.

 In 1979,  E. De Giorgi [3] conjectured that when a sequence of functionals,
 for instance, the one in (1.4) or in a more general form,
 converges in the sence of $\Gamma$-convergence to a limiting functional,
 the corresponding gradient flows will converge as well (maybe after an
 appropriate change of timescale). Also [4, p.216] and [5, p.507].

In [6], the author proved the De Giorgi's conjecture for a rather wide
kind of functionals.  Thus, a natural question is that under what conditions,
the functional sequence like (1.4) can be $\Gamma$-convergence.

A first result related to this question  was appeared in [7]. Because the integrands
in [7] have the same scale for the variables $x$ and $t$, the methods there can't
be applied to functionals (1.4) whose integrands are anisotropic in $x$ and $t$.

In this paper, we will cleverly combine the arguements in [8, 9, 10],  all
of which study the $\Gamma$-convergence of elliptic functionals like (1.3)
with the weak-topology of $W^{1,p}(\Omega , R^m)$,  to prove  that
the $\Gamma$-convergence holds for the functional (1.4)  under
assumption (1.5) and a periodic hypothesis (see (1.8) below). For this purpose, we
construct  functionals as follows.

 Let $Y={(0,1)}^n=\{0<y_i<1 , i=1,2, \cdots n\}$ , $kY={(0,k)}^n$, and $k_T=kY
 \times (0,T).$
\ For   $\lambda \in R^{mn} $ and a.e.  $ t \in R$, define
$$ \bar f(t,\lambda )=\inf_{k \in N}\inf \{|kY|^{-1}\int _{kY}
f(y,t,\lambda+D\phi(y,t))dy \colon \phi  \in V_p^0(k_T , m) \} , \eqno (1.6) $$
where and below $|E|\buildrel def \over = L^n(E)$ and $L^k$ is used to
denote the $k$-dimensional Lebesque measure.

 Obviously (1.5) implies that $\bar f(t,Du)$ is nonnegative and measurable, so we
can define the homogenized functional
$$ F(u, \Omega _T )=\int _{\Omega _T} \bar f(t, Du)dxdt , \ \ \ \   u \in V_p(\Omega _T ,
m).   \eqno (1.7) $$

 The  main result of this paper is the following theorem.

{\bf THEOREM 1.2.}   {\sl If hypotheses (1.4) and  (1.5) are satisfied,
and suppose
$$f(y, t, \lambda) \ \ is \ \ {\bar Y}-periodic \ \ on \ \
the \ \ first  \ \ variable \ \ y, \eqno (1.8)$$
 then for every
$T>0$ and every bounded open set $\Omega \subset R^n$ with $L^n( \partial \Omega )
=0$
$$ \Gamma (\tau ) \lim_{\varepsilon \to 0 }F^\varepsilon (u,\Omega _T) =
F(u,\Omega _T)   ,  \forall u \in V_p(\Omega _T ,m ) ,$$ }
 where $\tau $ is taken as the {\bf sw}-topology of $V_p(\Omega _T,m)$. (
 See the Def. 1.2 in [6] for the  {\bf sw}-topology.

The proof of this theorem will be given in section 4.

\

\

 \centerline { \bf 2.  PRELIMINARY LEMMAS }

We collect some
 properties of the $\Gamma$-limits in [1, 2] which are well-known but important
for  the coming arguements.

If the {\bf $ \limsup $ } in (1.2) is replaced by {\bf $\liminf$ }, the definition
1.1 is turned to the definition of {\bf low } $\Gamma $-limit. In this case, we denote
it by
$$ \lambda = \Gamma ^-(\tau )\lim_{h \to \infty }F^h(u) . $$
\noindent Similarly, we have {\bf upper} $\Gamma$-limit and denote it by
$\lambda =\Gamma ^+(\tau ) \lim_{h \to \infty }F^h(u) . $

Obviously, $\Gamma (\tau )\lim_{h\to \infty} F^h(u)$ exists if amd only if
$\Gamma ^+(\tau )\lim_{h\to \infty}F^h(u)=\Gamma^-(\tau )\lim_{h\to \infty } F^h(u) . $

{\bf LEMMA 2.1.}  {\sl $F^-(u) =\Gamma ^- (\tau )\lim_{h\to \infty}F^h(u) $ exists for
 every $ u \in X $, and $F^-(u)$ is lower semicontinuous in $ (X,\tau ) $.
If $ F(u)= \Gamma (\tau) \lim_{h \to \infty}F^h(u)$ exists for every $ u \in X$ ,
then F(u) is also lower semicontinuous  in $(X, \tau )$ .}

\vskip 0.4 cm

{\bf LEMMA 2.2.}  {\sl For each sequence $\{F^h\} $ of functionals in
$(X , \tau)$,
there exists a subsequence $F^{h_k}$ and $F^\infty $ from X to $\bar R$, such
 that $$F^\infty (u)=\Gamma (\tau )\lim_{k \to \infty }F^{h_k}(u) \quad  \forall
u \in X . $$ }

{\bf LEMMA 2.3.} {\sl Suppose that $\lambda =\Gamma (\tau ){\lim}_{\varepsilon
\to 0 }F^\varepsilon (u)$ and $\varepsilon _h \to 0 $ $(h \to \infty )$, then
$$\Gamma ^-(\tau ) \lim_{h \to \infty }F^{\varepsilon _h }(u)=\Gamma (\tau)\lim_{h \to \infty }F^{\varepsilon _h}(u) = \lambda $$ }

\vskip 0.4cm

{\bf LEMMA 2.4.}  {\sl Suppose that $ f \colon R\times R \to \bar R$, then
there exists a function $\delta \colon \varepsilon \to \delta (\varepsilon )$
such that $\varepsilon \to 0 $ implies $ \delta (\varepsilon ) \to 0 $ and
$$\limsup_{\varepsilon \to 0 }f (\delta (\varepsilon ),\varepsilon )
\le  \limsup_{\delta \to 0 } \limsup_{\varepsilon \to 0 }
f(\delta , \varepsilon )  \eqno (2.1) $$
Moreover, the opposite inequality for low limits and
the equality for limits hold true respectively. }

From now  on, we restrict ourselves to the sequence of functionals (1.4),
or more general
functionals :
$$F^\varepsilon (u,\Omega \times (a,b))=\int _a^b \int _{\Omega }f({x \over
\varepsilon },t,Du)dxdt, \ \ (\varepsilon \to 0^+).\eqno (2.2) $$

We will fix $T>0$ and allow $\Omega $ and $(a,b)$ to be arbitrary.
Let $S=R^n\times (0,T)$, $\beta _T$ be the $\sigma $ -ring generated by the set

$$\{ \Omega \times (a,b) \colon 0\le a<b \le T ,\Omega \subset R^n \ \ are \ \
bounded \ \  open \ \  sets \}. $$
\noindent Then $ (S,\beta _T ,L^{n+1})$ is a measure space. Let
$$V_{p,loc}=L^p([0,T], W_{loc}^{1,p}(R^n,R^m)) \eqno (2.3).$$

{\bf LEMMA 2.5.} {\sl Assume that (1.4), (1.5) and (1.8) are satisfied. Then for every
sequence $\varepsilon \to 0^+ $, there exist a subsequence $\varepsilon _h \to 0^+ $
 \ \ $(h \to 0 )$ and a family of $\sigma$-finite and $\sigma $-additive measures
$H(u,\Omega \times (a,b))$ on $\beta _T$, such that for every $ u \in V_{p,loc}$,
every finite interval $(a,b)$ and every bounded open set $\Omega \subset R^n$
with $L^n(\partial \Omega )=0$,
$$\Gamma (\tau )\lim _{h \to \infty }F^{\varepsilon _h}(u,\Omega \times (a,b))=
H(u, \Omega \times (a,b)) \eqno (2.4) $$ and
$$ 0 \le H(u,\Omega \times (a,b)) \le C\int _a^b \int _\Omega (1+|Du|^p)dxdt,
\eqno (2.5) $$
\noindent where $\tau $ is the {\bf sw}-topology of $V_p(\Omega \times (a,b)).$ }

{\bf \sl Proof.}   We follow the proof of in [9, Theorem 3.1]. D is used to
denote the algebra generated by all open cubes in $R^{n+1}$ with rational vertices
and E the class of all bounded open sets in $R^{n+1}$. Applying lemma 2.2 and a
diagonalization argument, we can find a sequence $ \varepsilon _h $ \ \ ($ h
\to \infty $) such that $\Gamma (\tau ) \lim_{h\to \infty} F^{\varepsilon _h}(u,Q) $ exists for
all $ Q \in D $, i.e
$$H^-(u,Q)=H^+(u,Q), \quad \forall Q \in D  ,$$
\noindent where $$H^-(u,Q)=\Gamma ^-( \tau )\lim_{h\to \infty} F^{\varepsilon _h}(u,Q)$$
and $$H^+(u,Q)=\Gamma ^+(\tau )\lim_{h\to \infty} F^{\varepsilon _h}(u,Q) . $$
\noindent In the same way as in [9, p.738-739], by lemma B in [6], we can prove that
$H^- $ is (finitely) super-additive and $H^+$ is sub-additive over D. For $e \in
E$, define
$$H(u,e)=\sup_{Q\subset \subset e}H^-(u,Q)=\sup_{Q\subset \subset e}H^+(u,Q) \quad Q\in D, $$
then $H(u,e)$ is an increasing, inner regular and finitely additive set function.
Therefor, the routine methods implies that (2.4) holds and $H(u,\Omega\times (a, b)
)$ can
be extended to a $\sigma $-finite and $\sigma $-additive measure on $ \beta _T$
(see  [11, Prop. 5.5 and Theorem 5.6 ]. From (2.4) and (1.5), the estimate (2.5)
follows immediately .

\

\

\centerline {\bf 3.  $\Gamma $-LIMITS OF LAYERED AFFINE FUNCTIONS }

 Throughout
this section,
suppose that (1.4), (1.5) and (1.8) are satisfied. $\tau$ is used to denote
the {\bf sw}-topology of $V_p(\Omega _T,m)$. For simiplicity, $V_p(\Omega _T)$
denotes the space $V_p(\Omega _T, m)$.
We intend to determine the $\Gamma $-limits of $F^\varepsilon (u,\Omega _T)$
for $u=\lambda (t)\cdot x+a(t) $ with $\lambda \in L^p([0,T],M(m\times n))$ and
$a \in L^p([0,T], R^m)$, where we define the norm on $M(m\times n)$,
the space of all real $m\times n$ matrices,  as the same as
on $R^{mn}$.

{\bf LEMMA 3.1.}   {\sl For each $u_{\lambda , a}=\lambda (t)\cdot x+a(t)$
with
$$\lambda \in L^p([0,T], M(m,n))  \ \ and \ \  a \in L^p([0,T], R^m), $$
\noindent there exists a sequence of functions $\{u^\varepsilon  \}\subset V_p(\Omega
_T)$ satifying
$$\{u^\varepsilon - u_{\lambda , a}\} \subset V_p^0(\Omega _T) \ \ and \ \
  u^\varepsilon \buildrel \tau \over \rightarrow
u_{\lambda ,a}  \ \ in \ \ V_p(\Omega_T) \ \ as \ \ \varepsilon \to 0^+ $$
such that
$$\lim_{\varepsilon \to 0^+}F^\varepsilon (u^\varepsilon , \Omega _T)=
\int _{\Omega _T}\bar f(t, \lambda )dxdt =
F(u_{\lambda , a}, \Omega _T), $$
\noindent where $\bar f(t,\lambda )$ is given by (1.6) and $F$ by (1.7).}

{\bf \sl Proof .}  Fix $ \delta \in (0,1)$, one can choose $ k \in N$
and $\phi ^\delta \in V_p^0(k_T ,m)$ (see (1.6)) such that
$$ \bar f(t,\lambda (t) ) \le |kY|^{-1}\int _{kY}f(y,t, \lambda +D\phi^\delta )
dy \le \bar f(t,\lambda (t)) + \delta . \eqno (3.1) $$
\noindent We use $E^{\ast}_\eta $ to denote the extension of $  \eta\bar{Y} $ on the
$\eta Y$-period, and  let
$$ \Omega^{\ast} _\eta = \{ e \in E^{\ast}_\eta , e \subset \Omega \} ,
\ \ \  E_{\eta}=\bigcup_{e \in E^\ast_\eta} e ,  \ \ \  \Omega_\eta =\bigcup
_{e \in \Omega^\ast_\eta} e ,$$
\noindent then $E_\eta = R^n $.  As $\Omega $ is bounded, $\Omega^\ast_\eta $ is
 a finite set for each
$\eta >0$ , and
$$\lim_{\eta \to 0^+ }L^n(\Omega \backslash \Omega _{\eta })= 0. \eqno (3.2) $$
\noindent For every $ t \in [0,T]$ , extend $ \phi ^\delta (y,t)$ such
that it is a $kY$-periodic function on the variable y, then define
$$ v^{\varepsilon , \delta }(x,t)= \cases {u_{\lambda , a}(x,t)+\varepsilon
\phi^\delta ({x \over \varepsilon }, t), &  $\Omega _{\varepsilon k}$ \cr
u_{\lambda ,a}(x,t) , & $ \Omega \backslash \Omega {\varepsilon k}$ \cr } .
\eqno (3.3) $$
\noindent It is easy to know that
$v^{\varepsilon , \delta } \in V_p(\Omega_T) , \ \  v^{\varepsilon , \delta }-u_{\lambda,a} \in V_p^0(\Omega _T).$
For each $D \in \Omega _{\varepsilon k}$, by the periodicity of
$$g(y,t)=f(y,t, \lambda (t)+D\phi ^\delta (y,t)), $$
\noindent we have
$$\eqalignno {\int _{D\times (0,T)} & f({x \over \varepsilon },t,Dv^{\varepsilon,\delta })dxdt
 =\int _0^T[\varepsilon ^n \int _{D/\varepsilon }f(y,t,\lambda (t)+D\phi ^\delta
(y,t)dy]dt  \cr
& = L^n(D)\int_0^Tdt |kY|^{-1}\int _{kY}f(y,t, \lambda (t)+D\phi^\delta )dy .& (3.4) \cr }$$
\noindent Summing up the both sides for all $ D \in \Omega _{\varepsilon k }$
and applying (3.1), we obtain that
$$\eqalignno {L^n(\Omega _{\varepsilon k})\int _0^T\bar f(t, \lambda (t))
dt & \le \int _0^Tdt \int _{\Omega _{\varepsilon k}}f({x \over \varepsilon },t,
Dv^{\varepsilon ,\delta })dx \cr
& \le L^n(\Omega _{\varepsilon k})\int _0^T(\bar f(t, \lambda (t))+ \delta )dt .  \cr  }$$
\noindent Thus, it follows from (1.5) and (3.3) that
$$\eqalignno { & L^n  (\Omega _{\varepsilon k}) \int _0^T \bar f(t, \lambda (t))
dt
  \le \int _{\Omega_T}f({x \over \varepsilon} ,t, Dv^{\varepsilon ,\delta })dx
dt \cr
& \le L^n(\Omega _{\varepsilon k})\int _0^T(\bar f(t, \lambda (t))+ \delta )dt + \int _0^T
dt\int _{\Omega \backslash \Omega _{\varepsilon k}}f({x \over \varepsilon },t,
\lambda (t))dx  \cr
& \le L^n(\Omega _{\varepsilon k})\int _0^T(\bar f(t, \lambda (t))+ \delta )dt +CL^n(\Omega
\backslash \Omega _{\varepsilon k})\int _0^T (1+|\lambda |^p)dt .
& (3.5) \cr }$$
\noindent By this estimate and (3.2), we see that
$$ \lim_{\delta \to 0^+}\lim_{\varepsilon \to 0^+}\int _{\Omega _T}f({x \over
\varepsilon},t,Dv^{\varepsilon ,\delta})dxdt=\int _{\Omega _T}\bar f(t, \lambda (t)
)dtdx . \eqno (3.6) $$
\noindent Moreover,  we have
$$ \|v^{\varepsilon ,\delta }-u_{\lambda ,a}\|_{L^p(\Omega _T)}^p
 = \varepsilon ^p \sum _{D \in \Omega _{\varepsilon k }}|D|\int _0^Tdt|kY|^{-1} \int
_{kY}|\phi ^\delta |^pdy  .  \eqno (3.7)   $$
\noindent Applying (3.6), (3.7), and lemma 2.4, one can find a sequence
$$\delta (\varepsilon ) \to 0^+ \ \ \ \ \  as \ \ \varepsilon \to 0^+ $$ such that
$\{u^\varepsilon = v^{\varepsilon, \delta(\varepsilon)} \colon  \varepsilon >0\}$
satisfy that
$$\{u^\varepsilon - u_{\lambda , a}\} \subset V_p^0(\Omega _T) , \ \
\lim _{\varepsilon \to 0^+} \|u^{\varepsilon }-u_{\lambda ,a}\|_{L^p(\Omega _T)}^p
=0$$
and
$$\lim_{\varepsilon \to 0^+} F^{\varepsilon}(u^{\varepsilon}, \Omega _T)
=F(u_{\lambda ,a}, \Omega _T).$$
On the other hand, the coercive condition in (1.5) and (3.5) imply that
$\{D u^\varepsilon \}$ is bounded in $L^p(\Omega _T, R^{mn})$. Thus, by
lemma B in [6], we obtain that
$  u^\varepsilon \buildrel \tau \over \rightarrow
u_{\lambda ,a}  .  $ This proves the desired result.

\vskip 0.4cm

{\bf LEMMA 3.2.}    {\sl Let $ u_{\lambda , a}(x,t) = \lambda (t)\cdot x+a(t) $
be the same as in lemma 3.1.
Then for each sequence $ u^\varepsilon \buildrel \tau \over
\rightarrow u_{\lambda , a} $
in $ V_p(\Omega _T) $  \ \ $( \varepsilon \to 0^+ )$,
$$  \liminf_{\varepsilon \to 0^+ }F^\varepsilon (u^\varepsilon , \Omega _T)
\ge F(u_{\lambda ,a} ,\Omega _T)=\int _{\Omega _T}\bar f(t,\lambda (t))dtdx .$$ }

{\bf \sl Proof .} {\bf  (1)}\ \    Firstly,  assume
 $ u^\varepsilon \buildrel \tau \over
\rightarrow u_{\lambda , a} $
and $ u^\varepsilon  - u_{\lambda , a}
\in V_p^0(\Omega _T) $. As $\Omega $ is bounded, we find an open cube D whose
sides are parallel to axes and whose center concides with the origin, such that
$ \bar \Omega \subset D $ . The side length of D is denoted by 2d, and let
$$k_\varepsilon = [ {2d \over \varepsilon } ]+3, \quad a_\varepsilon = [-{d \over
\varepsilon} ], $$
$$ x_\varepsilon = (a_\varepsilon , \cdots , a_\varepsilon ) \in R^n,
D_\varepsilon =\varepsilon (x_\varepsilon + k_\varepsilon Y ) ,$$
\noindent where $ [ \kappa ] $ denote the maximum integer not greater than $\kappa $.
It is not difficult to get
$$ D \subset D_\varepsilon , \quad  \lim_{\varepsilon \to 0^+}L^n
(D_\varepsilon ) = L^n(D) . \eqno (3.8) $$
\noindent Let $$ Q=D\backslash \bar \Omega ,  \quad  Q_T=Q \times (0,T) . \eqno
(3.9) $$
\noindent Applying lemma 3.1 to the open set $Q$, we can choose a sequence
$$v^\varepsilon \rightarrow u_{\lambda , a}  \ \  {\bf sw} \ \  in  \ \  V_p(Q_T), \quad
v^\varepsilon - u_{\lambda ,a} \in V_p^0(Q_T) $$
\noindent such that
$$\lim_{\varepsilon \to 0^+} F^\varepsilon (v^\varepsilon , Q_T) = \int _Q
\bar f(t,\lambda )dxdt . \eqno (3.10) $$
\noindent For fixed $ t \in  [0,T] , $  define
$$ \phi ^\varepsilon (x,t)= \cases {u^\varepsilon - u_{\lambda ,a} , & $ x \in
\bar \Omega$  \cr
v^\varepsilon - u_{\lambda ,a} , & $ x \in D\backslash \bar \Omega = Q $ \cr
0 , & $ x \in D_\varepsilon \backslash  D $  \cr } . \eqno (3.11) $$
By the periodicity of $ f(y,t, \lambda ),$
using the variable transformation, we obtain
$$\eqalignno{\int _{D_\varepsilon } & f({x \over \varepsilon } ,t, \lambda +D\phi ^\varepsilon
(x,t)) dx    =\varepsilon ^n\int _{x_\varepsilon +k_\varepsilon Y}f(y,t, \lambda +D_x\phi
^\varepsilon (\varepsilon y, t))dy  \cr
& = (k_\varepsilon \varepsilon)^n|k_\varepsilon Y|^{-1} \int_{k_\varepsilon Y} f(y,t, \lambda +D\psi
^\varepsilon (y,t))dy  ,
& (3.12) \cr }  $$
\noindent where $ \psi ^\varepsilon (y,t)=\varepsilon ^{-1} \phi ^\varepsilon
\bigl (\varepsilon (y+x_\varepsilon ),t\bigr ) $.
 Obviously, (3.11) gives us
 $$\psi ^\varepsilon \in V_p^0 \bigl ( (k_\varepsilon Y)\times (0, T)\bigr ).$$
Thus, we deduce, from
(3.12) and (1.6) yield that for each $t\in [0, T]$,
$$ |D_\varepsilon |^{-1}\int _{D_\varepsilon}f({x \over \varepsilon},t, \lambda +D\phi
^\varepsilon (x) )dx =|k_\varepsilon Y|^{-1} \int _{K_\varepsilon Y}f(y,t, \lambda +D\psi ^\varepsilon
(y,t))dy \ge \bar f(t,\lambda )  .  $$
Therefore
$$ \int_0^T  \int _{D_\varepsilon}f({x \over \varepsilon},t, \lambda +D\phi
^\varepsilon)dxdt \ge L^n(D_\varepsilon) \int _0^T\bar f(t,\lambda )dt . $$
On the other hand , by (1.5) and (3.8), we have
$$ \liminf_{\varepsilon \to 0^+}\int _0^Tdt \int _{D_\varepsilon}f({x
\over \varepsilon} , t, \lambda +D\phi ^\varepsilon )dxdt =\liminf_{
\varepsilon \to 0^+} \int _0^T dt \int _Df({x \over \varepsilon },t, \lambda +
D\phi^\varepsilon )dx .$$
This yields  $$ \liminf_{\varepsilon \to 0^+}F^\varepsilon\bigl ( u_{\lambda ,a}+\phi
^\varepsilon , D\times (0,T)\bigr ) \ge \int _0^T\int_D \bar f(t,\lambda )dxdt . $$
Combing this estimate, (3.9), (3.10) with (3.11), we have
$$ \eqalignno{\liminf_{\varepsilon \to 0^+}F^\varepsilon (u^\varepsilon ,
\Omega _T) &  = \liminf_{\varepsilon \to 0^+}[F^\varepsilon (u_{\lambda ,a}
+\phi ^\varepsilon , D\times (0,T)) - F^\varepsilon (v^\varepsilon ,
Q\times (0,T)) ] \cr
& \ge\int_0^T \int_D\bar f(t,\lambda )dxdt - \int_0^T\int _Q\bar f(t, \lambda )dxdt \cr
& = \int _{\Omega _T}\bar f(t, \lambda )dxdt . \cr } $$
{\bf (2)}\ \   In order to remove the restriction
$ u^\varepsilon -  u_{\lambda ,a} \in V_p^0(\Omega _T) ,  $
it is sufficient to apply the  {\it De Giorgi's } arguements and the
result of the case (1). See [11], or [8, p.197] for the details.

\vskip 0.4cm

{\bf DEFINITION 3.3.}   {\sl Let $ \{ \Omega _i \colon i= 1 ,2, \cdots h \}$ be a finite
partition of $ \Omega $ into open sets (except for a set of measure zero ),
$\lambda _i \in L^p([0,T],M(m,n)) , \ \  a_i \in L^p([0,T], R^m) . $  We call the function
$$ W(x,t) = \cases {\lambda _i (t) \cdot x +a_i(t) , & $x\in \Omega _i $  \cr
0 , & $x \in \Omega \backslash \cup _{i=1}^h \Omega _i $ \cr } $$
a $L^p$ -layered affine function on $ \Omega _T$. }

\vskip 0.4cm

To sum up lemmas 3.1 and 3.2 (observing that $\Omega $ maybe arbitrary there ),
lemmas 2.5 and 2.3, we obtain the following theorem.

{\bf THEOREM 3.4.}    {\sl Suppose that $\Omega $ is a bounded open set in $R^n$
with $L^n(\partial \Omega )$ $=0$, $H(u,\Omega _T) $ is given by lemma 2.5, then
$$\Gamma (\tau ) \lim_{\varepsilon \to 0^+} F^\varepsilon (w, \Omega _T)=\int
_{\Omega _T} \bar f(t,Dw)dxdt = H(w, \Omega _T) $$
for any $ w$, a  $L^p$ -layered affine function  on $\Omega _T$. }

\

\

\centerline {\bf 4 . A PROOF OF THEOREM 1.2 }

In this section, we suppose that all  the hypotheses of Theorem 1.2 are satisfied.
Applying  the same arguement as  in [9, Section 5], we can
prove that for almost
$ t \in [0,T] $, $ \bar f(t, \lambda )$ is convex if $n=2$; and convex with
respect to each column vector if $ n >2$. This implies that

{\bf LEMMA 4.1.}  {\sl For a.e. $ t \in [0,T] $, $ \bar f(t, \lambda )$ is
continuous in $M(m,n)$ .}

\vskip 0.4cm

{\bf LEMMA 4.2.}   {\sl Suppose $v \in V_p(\Omega _T) \ \ \ (1<p<\infty )$, then
there exists a sequence of $ L^p$-layered affine functions :}
$$ v^k(x,t) = \cases {\lambda_i^k(t)\cdot x +a_i^k(t) , & $x \in \Omega_i$ \cr
0 , & $x \in \Omega \backslash \cup_{i=1}^{h_k}\Omega_i $ \cr } $$
such that
$\|v-v^k \|_{V_p(\Omega_T)} \longrightarrow 0 \ \ \ as \ \ \  k \to \infty $.

{\bf \sl Proof .}   {\bf (1)} \ \ Suppose $ 1<p \le 2 $. Fix $v \in V_p(\Omega_T) $. For
any $\varepsilon >0$ we can choose $ u \in V_2(\Omega_T)$ such that
$$ \| u-v\|_{V_p(\Omega_T)} < \varepsilon . \eqno (4.1) $$
Because $ H \buildrel \rm def \over = W^{1,2}(\Omega)$ is a {\sl Hilbert }
space, one can assume that $ \{\psi_l\}_{l=1}^\infty $ is its complete orthonormal
basis. Let
$$C_l(t)= <u(t,\cdot ) , \psi _l >_H , $$
then $ C_l(t) \in L^2[0,T] $, and for a.e. $ t \in [0,T] $
$$I_k(t)=\|u -\sum_{l=1}^kC_l\psi_l(x) \|_H \longrightarrow 0 \ \ \ ( k \to \infty ).$$
Thus the domainnated convergence theorem implies that for some integer k
$$ \| u -\sum_{l=1}^k C_l\psi_l \|_{V_p(\Omega_T)} \le \varepsilon . \eqno (4.2) $$
It is well known that there exist piecewise affine functions $ \omega_l(x) $
in $\Omega $ such that
$$\max_{1\le l \le k}\|\psi_l - \omega_l \|_{W^{1,p}(\Omega)} \le \varepsilon (1+
 \sum_{l=1}^k\|C_l\|_{L^p(\Omega)}^{-1}) .$$
Let  $$v^\varepsilon (x,t) =\sum_{l=1}^k C_l(t) \omega _l (x) , $$
then
$$\| v^\varepsilon - \sum_{l=1}^k C_l\psi _l \|_{V_p(\Omega_T)} \le C(p)
\varepsilon . \eqno (4.3) $$
Combing (4.1), (4.2) with (4.3), we get

$$ \| v-v^\varepsilon \|_{V_p(\Omega_T)} \le C(m,n,p) \varepsilon .$$
Observing that each $v^\varepsilon $ can be written as a layered function on
$\Omega_T$, we have completed the proof.

{\bf (2)} \ \  Suppose  $2<p< \infty $.   Applying {\sl Sobolev }embedding theorem
we can find an integer
$ k , \ \  {k-1 \over n } \ge {1 \over 2} -{1 \over p},$
such that
$$ H_1 \buildrel \rm def \over = W^{k,2}(\Omega) \hookrightarrow W^{1,p}(\Omega) .$$
Given $ v\in V_p(\Omega_T) $. For $ \varepsilon >0$ , one can find
$ u \in L^p([0,T],W^{k,2}(\Omega))$
such that
$$ \| u-v \|_{V_p(\Omega_T)} < \varepsilon  .  \eqno (4.4)  $$
Let $ \{ \psi _l\}_{l=1}^\infty $ be the complete orthonormal basis of the
{\sl Hilbert } space $ H_1$, then
$$C_l(t)\buildrel def \over  = <u(\cdot,t) , \psi_l>_{H_1} \in L^p[0,T] .$$
The remaining part is entirely the same as the case (1).

\vskip 0.4cm

Now we are in the position to {\bf prove  Theorem 1.5. }  We will use the idea of
[9, p.750-751]. For $ u \in V_p(\Omega_T) $, we can extend $u$ such
that $ u \in V_{p,loc} $ (recall (2.3) ). From lemma 4.2, choose a sequence of
$L^p$ -layered functions $\omega^k(x,t)$, such that
$$ \| u- \omega^k \|_{V_p(\Omega_T)} \longrightarrow 0  \ \ (k \to \infty ) .
\eqno (4.5) $$
By taking a subsequence, one can assume that $ D\omega^k \to Du $ almost everywhere
on $\Omega_T$ and
$$ \bar f(t ,D\omega^k) \longrightarrow \bar f(t, Du)  \ \ a.e \ \ in \ \ \Omega_T$$
by the virtue of the continuity of $\bar f(t,\cdot )$ (see Lemma 4.1 ).

We deduce, from the absolute continuity of $ \int |Du|^pdxdt $, {\sl Egoroff}
theorem, theorem 3.4 and  inequality (2.5), that
$$ \liminf_{k \to \infty}\int_{\Omega_T}\bar f(t, D\omega^k)dxdt
\le \int_{\Omega_T}\bar f(t, Du )dxdt . $$
Therfore, by the semi-continuity of $H(u,\Omega_T) $ (see lemma 2.1 ),
$$\eqalignno {H(u, \Omega_T) &  \le \liminf_{k \to \infty}H(\omega^k ,
\Omega_T)  \cr
& =\liminf_{k \to \infty}\int_{\Omega_T}\bar f(t,D\omega^k)dxdt \cr
& \le \int_{\Omega_T}\bar f(t,Du)dxdt . &(4.6) }$$
On the other hand,  according to lemma 2.5 and {\sl Lebesgue-Nikodym} theorem
(see \S 3 of Ch.3 in [12] ), we have
$$ H(u,\Omega_T)=\int_{\Omega_T}h(x,t)dxdt  \eqno (4.7) $$
for some $h \in L_{loc}^1(R^n \times (0,T) )$ and all $ \Omega_T=\Omega \times
(0,T) $.  By approximation argument, one can easily prove that for a.e
$(x,t) \in \Omega_T $, there exists $ r_k \to 0^+$ such that
$${u(x+r_k(y-x),t)-u(x,t) \over r_k} \buildrel \tau \over \rightarrow Du(x,t)\cdot (y-x) \ \
 \ \ in \ \ V_p(B(x,1)\times (0,T)) .  \eqno (4.8) $$
Since
$$ |\int_0^Th(x,t)dt - \int_0^T |B(x,r_k)|^{-1}\int_{B(x,r_k)}h(y,t)dydt |$$
$$\le |B(x,r_k)|^{-1}\int_{B(x,r_k)}
|\int_0^T[h(x,t)-h(y,t)]dt|dy $$
and
$$ \int_0^Th(y,t)dt \in L_{loc}^1(R^n) ,$$
so
$$ \int_0^Th(x,t)dt=\lim_{k \to \infty}\int_0^T |B(x,r_k)|^{-1}\int_{B(x,r_k)}h(y,t)dydt
\ \  for \ \ a.  e. \ \  x \in \Omega .  \eqno (4.9) $$
Fix $ k \ \ and \ \ x $, set
$$r=r_k,   \ \  B_r=B(x,r),   \  \   B_{r,T}=B_r \times (0,T) .$$
By (4.7) and lemma 2.5, we can find a sequence
$$u^h \buildrel \tau \over \rightarrow u   \ \  in \ \  V_p(B_{r,T})  (h\to 0)$$
such that
$$ \eqalignno {\int_0^T  & |B_r|^{-1}\int_{B_r}h(y,t)dydt =|B_r|^{-1}H(u,B_{r,T})  \cr
& =\lim_{h \to \infty}\int_0^T |B_r|^{-1}\int_{B_r}f({y+\varepsilon_hk_h \over \varepsilon_h} ,
t,Du^h)dydt, \ \  \Bigl (k_h \buildrel \rm def \over = [{x(r-1) \over \varepsilon_h}]
\Bigr )  \cr
& \ge \liminf_{h \to \infty}\int_0^T |B_{{r\over 2}}|^{-1}\int_{B_{r \over 2}}f({y+x(r-1) \over
\varepsilon_h },t, Du^h(y+a_h ,t))dydt  \cr
&  \ \ \ \ \Bigl ( note \ \  that \ \  a_h \buildrel \rm def \over =
x(r-1)-\varepsilon_hk_h \to 0^+ \Bigr )  \cr
& =\liminf_{h \to \infty}\int_0^T |B_{{1\over 2}}|^{-1}\int_{B_{1 \over 2}}f
\Bigl ( {ry \over \varepsilon
_h} ,t, D\bigl ( u^h(x+r(y-x)+a_h ,t) \cr
& \ \ \ \ \ \ \ \ \ \ \ \ -u(x,t) \big ) r^{-1}\Bigr )dydt .
& (4.10)  \cr }$$
Let $u_{r,x}(y,t)=r^{-1}[u(x+r(y-x), t)-u(x,t)]$.  Obviously
$$r^{-1}[u^h(x+r(y-x)+a_h ,t) -u(x,t) ] \buildrel \tau \over \rightarrow r^{-1}u_{r,x}
\ \  in \ \ V_p(B_{{1 \over 2},T}) \ \  as \ \  h \to \infty . $$  Let
$$\delta_h=r^{-1} \varepsilon_h  ,  \ \ \ \   a=|B_{1 \over 2}|^{-1},
\ \ \ \ \
F^-(u,Q) = \Gamma^-(\tau)\lim_{h\to \infty}F^{\delta_h}(u,Q) . $$
By (4.10), lemmas 2.1 and 2.3, (4.8) and theorem 3.4 in that order,
we deduce that
$$\eqalignno { \lim_{k \to \infty}\int_0^T |B(x,r_k)|^{-1}\int_{B(x,r_k)}h(y,t)dydt & \ge
a \cdot \liminf_{k \to \infty }F^-(u_{r_k,x}, B_{{1 \over
2},T} )  \cr
& \ge a \cdot  F^-(Du(x,t)(y-x) ,B_{{1 \over 2},T} )   \cr
& =a\int_0^T \int_{B(x,{1 \over 2})}\bar f(t,Du(x,t) dt dy    \cr
& =\int_0^T \bar f(t,Du(x,t))dt  . \cr } $$
Combing this estimate with (4.9), we obtain
$$ \int_{\Omega_T} h(x,t)dxdt \ge \int_{\Omega_T}\bar f(x,Du)dxdt  , $$
which together with (4.7) implies the opposite inequality of (4.6). Hence
$$ H(u,\Omega_T) = \int_{\Omega_T} \bar f(t,Du) dxdt  \ \  \forall  u \in V_p
(\Omega_T) . $$
Observing that
$$ F(u , \Omega_T) = \int_{\Omega_T}\bar f(t,Du)dxdt $$
is independent of $\{ \varepsilon_h\}$, we have completed the proof of Theorem
1.2 by lemma 2.5.

\

\

\

\

\centerline{\bf REFERENCES }

\item{[1]} H. Attouch, Variational convergence for functions and operators,
 Appl. Math. Series, Pitman, 1984 .

\item{[2]} G. Dal Maso, An introduction to $\Gamma$-convergence, Boston Basel:
Birkh\"auser, 1993.

\item{[3]} E. De Giorgi, New problems in $\Gamma$-convergence and G-convergence
in Free Boundary Problems, Proc. of seminar held in Pavia, September-October
1979, Ist. Naz. Alt. Mat. Francesco Severi, Vol. II, 183-194, Rome, 1980.

\item{[4]}L. Bronsard and R.V. Kohn, Motion by mean curvature as the singular limit
of Ginzburg-Landau Dynamics, {\sl J. Differ. Eqn.} {\bf 90} (1991), 211-237.

\item{[5]} N. C. Owen, J. Rubinstein and P. Sternberg, Minimizers and gradient
flows for singularly perturbed bi-stable potentials with a Dirichlet
condition, {\sl Proc. R. Soc. Lond Ser.A} {\bf429} (1990), 505-532.

\item{[6]} H.Y. Jian, A relation between $\Gamma$-convergence of functionals and
their associated gradient flows, Science in China, Ser. A,
42(2)(1999), 133-139.

\item {[7]} H.Y. Jian, Homogenization problems of parabolic minima, {\sl
Acta Math. Appl. Sinica} {\bf 12 } (1996), 318-327.

\item{[8]} S. M\"uller, Homogenization of noncovex integral functionals and
elastic materials, {\sl Arch. Rat. Mech. Anal.} {\bf 99} (1987), 187-212.

\item{[9]} E. Weinan, A class of homogenization problems in the calculus of variations,
{\sl Commu. pure Appl. Math.} {\bf 44} (1991), 733-759.

\item{[10]} H.Y. Jian, $\Gamma$-convergence of noncoercive functionals in
vector-valued space $W^{1,p}(\Omega )$ (in chinese), {\sl Science in China Ser. A}
{\bf 24(3) } (1994), 233-240.

 \item{[11]} E. De Giorgi  and G. Letta,  Une notion general de convergence
faible pour des fontions crotssantes d'ensemble, {\sl Ann  Scuola Norm. Sup. Pisa
CI  Sci.}  {\bf 4} (1977), 61 -99.

\item{[12]} K. Yosida,  Functional analysis, 5th Ed, Springer-Verlag,
1975.

\end